\documentclass{article}
\usepackage{amssymb}
\newcommand{\dem}{\noindent{\bf Proof. }}
\newcommand{\fg}{{\mathfrak g}}

\newcommand{\ad}{\mbox{\rm ad}}

\newcommand{\nuc}{\mbox{\rm ker}}
\newcommand{\id}{\mbox{\rm id}_\fg}

\newcommand{\f}{\mathfrak}
\newcommand{\qed}{\hfill $\square$}
\newcommand{\tR}{\tilde{R}}

\newtheorem{theor}{Theorem}[section]
\newtheorem{prop}[theor]{Proposition}
\newtheorem{coro}[theor]{Corollary}
\newtheorem{lemma}[theor]{Lemma}
\newtheorem{rema}{Remark}[section]
\newtheorem{defi}{Definition}[section]
\newtheorem{example}{Example}[section]
\newtheorem{examples}{Examples}[section]

\begin{document}

\title{Pseudo-K\" ahler Lie algebras with abelian complex structures}

\author{Ignacio Bajo$^{(1)}$ and Esperanza Sanmart\'{\i}n$^{(2)}$}

\maketitle

\begin{enumerate}
\item[(1)] Depto. Matem\'atica Aplicada II, E.I. Telecomunicaci\'on, Campus Marcosende, 36310 Vigo, Spain.

 ibajo@dma.uvigo.es
\item[(2)]  Depto. Matem\'aticas, Facultad de CC.EE., Campus Marcosende, 36310 Vigo, Spain.

esanmart@uvigo.es
\end{enumerate}

\begin{center}{\bf Abstract}\end{center}
{\small We study Lie algebras endowed with an abelian complex structure which admit a symplectic form compatible with the complex structure. We prove that each of those Lie algebras is completely determined by a pair $(U,H)$ where $U$ is a complex commutative associative algebra and $H$ is a sesquilinear hermitian form on $U$ which verifies certain compatibility conditions with respect to the associative product on $U$. The Riemannian and Ricci curvatures of the associated pseudo-K\"ahler metric are studied and a characterization of those Lie algebras which are Einstein but not Ricci flat is given. It is seen that all  pseudo-K\"ahler Lie algebras can be inductively described by a certain method of double extensions applied to the associated complex asssociative commutative algebras.}

\section{Introduction}

An abelian complex structure on a Lie algebra $\fg$ is a linear map $GL(\fg)$ such that $J^2=-\mbox{id}_\fg$ and $[Jx,Jy]=[x,y]$ for all $x,y\in\fg$. It is obvious that the Nijenhuis tensor for such a $J$ vanishes and, therefore, the structure is integrable. The name {\it abelian} comes from the fact that the eigenspaces of $J$ in the complexification $\fg^{\mathbb C}$ of $\fg$ are abelian Lie algebras. Note that they are a particular case of the so called nilpotent complex structures \cite{CFU}.

A pseudo-K\"ahler Lie algebra is a symplectic Lie algebra $(\fg,\omega)$ with a complex structure $J$ which is skew-symmetric with respect to $\omega$. Such a Lie algebra is equipped with the pseudo-Riemannian metric $g$ defined by $g(x,y)=\omega (Jx,y)$ for all $x,y\in\fg$ for which the complex structure is parallel. The aim of this paper is to study the structure and main properties of Lie algebras which are pseudo-K\"ahler for an abelian complex structure.

If an algebra $\fg$ admits an abelian complex structure it must be 2-step solvable. One can easily verify that Lie algebras of dimension 2 have an abelian complex structure (actually, all complex structures on these algebras are abelian) and that every non-degenerate skew-symmetric bilinear form $\omega$ is a symplectic form compatible with the complex structure. So, 2-dimensional Lie algebras are the first examples of pseudo-K\"ahler Lie algebras with abelian complex structure. Recall that the unique non abelian 2-dimensional Lie algebra is the algebra $\f{aff}({\mathbb R})$ of affine motions of ${\mathbb R}$. An important role in the sequel will be played by the Lie algebras $\f{aff}({\mathcal A})$ constructed as the tensor product $\f{aff}({\mathcal A})=\f{aff}({\mathbb R})\otimes_{\mathbb R}{\mathcal A}$ of $\f{aff}({\mathbb R})$ with an associative commutative algebra ${\mathcal A}$. A well-known example of these algebras is the underlying Lie algebra of the Kodaira-Thurston manifold, which is the algebra $\f{aff}({\mathcal A})$ when ${\mathcal A}$ is the 2-dimensiona nilpotent power algebra. All the algebras $\f{aff}({\mathcal A})$ carry abelian complex structures but in order to admit a compatible symplectic form the algebra ${\mathcal A}$ must verify further conditions. The existence of a symmetric bilinear $B$ form on  ${\mathcal A}$ such that the pair $({\mathcal A}, B)$ is a symmetric algebra \cite{curtis} guarantees that the corresponding Lie algebra $\f{aff}({\mathcal A})$ admits a pseudo-K\"ahler structure for an abelian complex structure; however, these examples do not exhaust the family of pseudo-K\"ahler Lie algebras with abelian complex structure.

In this paper we will show that every pseudo-K\"ahler Lie algebra with abelian complex structure is completely characterized by a pair $(U,H)$ where $U$ is a complex associative commutative algebra and $H$ a sesquilinear hermitian form on $U$ which verifies certain compatibility conditions. This characterization let us calculate some nice formulas for the Riemannian and the Ricci curvatures of the pseudo-K\"ahler metric and lead us to a complete description those algebras in our family for which the metric is Einstein but not Ricci flat. Moreover, in the last section we prove that all the pseudo-K\"ahler Lie algebras with abelian complex structure can be obtained by sucessive application of a method of double extension, which consist of a central extension and a generlized semi-direct product, on the associated pairs $(U,H)$ starting from an algebra $\f{aff}({\mathcal A})$.

\section{Preliminaries}

We first recall some basic definitions \cite{al-m-t}, \cite{Andrada1}, \cite{goze}, \cite{helm}. All the
algebras  in the paper are considered real or complex finite dimensional algebras.

\begin{defi} \label{def1} {\em A {\it complex structure} on a real Lie algebra $\fg$ is a linear map $J\in {\mathfrak{gl}}(\fg)$ 
 such that $J^2=-\id$ and
$$[Jx,Jy]=[x,y]+J[Jx,y]+J[x,Jy]\quad \mbox{for all }x,y\in\fg. $$

Two Lie algebras endowed with complex structures $(\fg_1,J_1),(\fg_2,J_2)$ are said to be {\it holomorphically
equivalent} if there is an isomorphism $\psi:\fg_1\to\fg_2$ such that $J_1=\psi^{-1}\circ J_2\circ \psi$.

We will say that a complex structure $J$ is {\it abelian} if $[Jx,Jy]=[x,y]$ holds for all $x,y\in\fg.$ }
\end{defi}

\begin{defi} {\em Let $\fg$ be a complex or a real Lie algebra. A {\it product structure} on $\fg$ is a linear map $K\in {\mathfrak{gl}}(\fg)$ 
 such that $K^2=\id$, $K\neq\id$ and
$$[x,y]=K[Kx,y]+K[x,Ky]-[Kx,Ky]\quad \mbox{for all }x,y\in\fg. $$
The product structure $K$ provides a decomposition of the vector space $\fg$ as the direct sum of 
the eigenspaces $\fg_+=\nuc (K-\id),\fg_-=\nuc (K+\id)$. Actually, $\fg_+$ and $\fg_-$ are subalgebras of 
$\fg$. When both subalgebras $\fg_\pm$ have the same dimension, the product structure is said to be a 
{\it paracomplex structure} on $\fg$.

We will say that a product structure $K$ is {\it abelian} if $[Kx,Ky]=-[x,y]$ holds for all $x,y\in\fg.$ In such 
case one easily verifies that the subalgebras $\fg_\pm$ are actually abelian.}
\end{defi}

\begin{rema}{\em One can also use Definition \ref{def1} to define the notion of complex structure for complex Lie algebras, considering $J$ to be ${\mathbb C}$-linear. However, this is not very interesting since for a complex Lie algebra $\fg$ there is one-to-one correspondence between such ${\mathbb C}$-linear complex structures and product structures on $\fg$. Actually, a straightforward computation shows that a ${\mathbb C}$-linear complex mapping $J$ on a complex Lie algebra $\fg$ is a complex structure if and only if the map $K=iJ$ is a product structure, where $i$ stands for the imaginary unit. Further, $J$ is abelian if and only if $K$ is so.}\end{rema}

\medskip

The following result is very well-known (see, for instance, \cite{Andrada1}):

\begin{lemma} Every Lie algebra admitting an abelian complex structure is 2-step solvable. \end{lemma}

\begin{examples}{\em $\,$ Some interesting examples of Lie algebras admitting abelian complex structures are:

\begin{enumerate}
\item {\it The Lie algebra $\mathfrak{aff}({\cal A})$ of a commutative associative algebra}. 

Let ${\cal A}$ be a commutative associative algebra. 
The vector space ${\cal A}\oplus {\cal A}$ with the product defined by
$[(a,b),(a',b')]:= (0,ab'-a'b)$ for all $a,a',b,b' \in {\cal A},$
is a Lie algebra denoted by $\f{aff}({\cal A})$ \cite{BD1}. Note that, actually, $\f{aff}({\cal A})=\f{aff}({\mathbb R})\otimes {\cal A}$ with the bracket $[x\otimes a,y\otimes a']=[x,y]\otimes aa'$. It is clear that the linear map  
$J$ on $\f{aff}({\cal A})$ defined by: $J(a,b):= (-b,a)$, for $a,b \in {\cal A},$ is an abelian
 complex structure on $\f{aff}({\cal A}).$ 

\item {\it The underlying real algebra of a complex Lie algebra with an abelian product structure}

When $(\fg, K)$ is a complex Lie algebra with a ${\mathbb C}$-linear abelian product structure, the underlying real Lie algebra $\fg_{\mathbb R}$ is naturally endowed with the abelian complex structure defined by the ${\mathbb R}$-linear map $J=iK$. 
\end{enumerate}
}\end{examples}

\begin{defi} {\em We say that a Lie algebra $\fg$ admits a {\it symplectic structure} if it admits a non-degenerate scalar
 2-cocycle $\omega$. The pair $(\fg,\omega)$ is said to be a {\it symplectic Lie algebra}. Two symplectic Lie algebras  $(\fg_1,\omega_1)$,  $(\fg_2,\omega_2)$ are said to be {\it symplectomorphic} if there exists a Lie algebras homomorphism $\varphi:\fg_1\to\fg_2$ such that $\omega_1(x,y)=\omega_2(\varphi x,\varphi y)$ for all $x,y\in\fg_1$.}
\end{defi}

\begin{rema} {\em A symplectic Lie algebra $(\fg,\omega)$ is naturally endowed with a
structure of left-symmetric algebra compatible with the Lie structure; this is to say, the product defined by $\omega (x\cdot y,z)=-\omega (y,[x,z])$ for $x,y,z\in\fg$
verifies $[x,y]=x\cdot y-y\cdot x$ and
$$[x,y]\cdot z=x\cdot (y\cdot z)-y\cdot (x\cdot z)\quad \mbox{for all }x,y,z\in\fg. $$
Geometrically, this means that each Lie group with Lie algebra $\fg$ can be equiped with the torsion free flat left-invariant 
connection defined by $\nabla^\omega_xy=x\cdot y$ for $x,y\in\fg$. }
\end{rema}

\begin{defi} {\em Let $(\fg,\omega)$ be a symplectic Lie algebra and $J$ a complex structure on $\fg$ such that 
 $\omega(Jx,Jy)=\omega(x,y)$ for all $x,y\in\fg$. The triple $(\fg,\omega,J)$ will be called a {\it
pseudo-K\"ahler Lie algebra}. If the complex structure is abelian we will say that the algebra $\fg$ is 
endowed with an {\it abelian pseudo-K\"ahler structure}. 
The pseudo-Riemannian metric $g$ defined on $\fg$ by $g(x,y)=\omega (Jx,y)$ will be called the {\it pseudo-K\"ahler metric} of $(\fg,\omega, J)$.
}\end{defi}

The following lemma will be useful in the next sections.

\begin{lemma} \label{commutes} Let $(\fg,\omega,J)$ be a pseudo-K\"ahler Lie algebra where $J$ is abelian. The left-symmetric product defined
by $\omega$ on $\fg$ verifies $J(x\cdot y)=(Jx)\cdot y$ for all $x,y\in\fg$.
\end{lemma} 
\dem If $x,y,z\in\fg$ then we have
$$\omega (J(x\cdot y),z)=-\omega (x\cdot y,Jz)=\omega(y,[x,Jz])=-\omega(y,[Jx,z])=\omega((Jx)\cdot y,z),$$
which proves the result.\qed

\begin{rema} {\em One may think that for every associative and commutative algebra ${\cal A}$ the corresponding Lie algebra $\f{aff}({\cal A})$ always admits an abelian pseudo-K\"ahler structure. However the result is not true and one may find associative commutative algebras ${\cal A}$ such that $\f{aff}({\cal A})$ does not even admit a symplectic form. For instance, if ${\cal A}$ is the 3-dimensional associative algebra spanned by $a_1,a_2,a_3$ with non-trivial products $a_1a_1=a_2a_2=a_3$, then the algebra $\f{aff}({\cal A})$ is spanned by $E_1,\dots ,E_6$ where $E_j=(a_j,0)$ and $E_{3+j}=(0,a_j)$ for $1\le j\le 3$ and their only non-trival brackets are $[E_1,E_4]=[E_2,E_5]=E_6$. If $\omega$ is a 2-cocycle on $\f{aff}({\cal A})$, one has
\begin{eqnarray*}
& & \omega (E_6,E_j)=\omega ([E_1,E_4],E_j)=-\omega ([E_4,E_j],E_1)-\omega ([E_j,E_1],E_4)\\
& & \omega (E_6,E_j)=\omega ([E_2,E_5],E_j)=-\omega ([E_5,E_j],E_2)-\omega ([E_j,E_2],E_5).
\end{eqnarray*}
The first identity shows that $\omega (E_6, E_2)=\omega (E_6, E_3)=\omega (E_6, E_5)=0$ and the second one implies $\omega (E_6, E_1)=\omega (E_6, E_4)=0$. Therefore $\omega(E_6,x)=0$ for all $x\in \f{aff}({\cal A})$, which proves that the 2-cocycle cannot be non-degenerate.
}\end{rema}

We recall the following definition  \cite{curtis}:
\begin{defi} \label{def_symm} {\em
An associative algebra ${\cal A}$ over ${\mathbb K}$ is said to be a {\it symmetric 
algebra} if ${\cal A}$ is equipped with a non-degenerate symmetric bilinear form $B: {\cal A}\times {\cal A} \rightarrow {\mathbb K}$  that satisfies
 $B(ab,c)= B(a,bc)$ for $a,b,c\in {\cal A}$. }
\end{defi}

It should be noticed that, according to the results given in \cite{borde}, one can construct a symmetric associative commutative algebra by applying the method of T$^*$-extension to an arbitrary associative commutative algebra.

\bigskip

\begin{examples}\label{firstexam} {\em $\,$ The following examples of pseudo-K\"ahler Lie algebras with abelian complex structures show that the class of such algebras is a wide one.

\begin{enumerate}
\item {\it The Lie algebra $\mathfrak{aff}({\cal A})$ of a symmetric associative commutative algebra}. 

If $({\cal A},B)$ is a commutative associative symmetric algebra then the algebra $\f{aff}({\cal A})$ with the abelian complex structure defined above and the symplectic form given by
$$\omega((a,b),(a',b')):= B(a,b')-B(b,a')\, ,\quad a,a',b,b' \in {\mathcal A}$$
is pseudo-K\"ahler. In this case we will say that $\f{aff}({\cal A})$ is equipped with the {\it standard pseudo-K\"ahler structure}.

A simple calculation shows that the left-symmetric product defined by $\omega$ on $\mathfrak{aff}({\cal A})$ is given by
$$(a,b)\cdot (a',b')=-(aa',ba')\, ,\quad a,b,a',b'\in {\mathcal A}\, ,$$
which is, actually, associative.

\item {\it The underlying real algebra of a complex Lie algebra with an abelian para-K\"ahler structure}

If $(\fg, K,\omega)$ is a complex Lie algebra with an abelian  para-K\"ahler structure, then the underlying real Lie algebra $\fg_{\mathbb R}$ is pseudo-K\"ahler if one considers the complex structure $J=iK$ and the symplectic form $\omega_{\mathbb R}=\mbox{Re}(\omega)$ where $\mbox{Re}(\omega)$ stands for the real part of $\omega$.

It should be noticed that the real dimension of each of those algebras is $4k$ for some $k\in {\mathbb N}$.

\item {\it An example with non associative left-symmetric product}

Let us consider the 6-dimensional real Lie algebra $\fg$ linearly spanned by $\{E_1,\dots,E_6\}$ with the non-trivial brackets
$$\begin{array}{lll}
[E_1,E_2]=-E_3-\frac{1}{2}E_6\, , & [E_1,E_4]=4E_2\, , &
{[E_1,E_5]}=\frac{1}{2}E_3+E_6\, ,\vspace{0.15cm}\\
{[E_2, E_4]}=\frac{1}{2}E_3+E_6\, , & [E_4, E_5]=-E_3-\frac{1}{2}E_6\, . &
\end{array}$$
Notice that $\fg$ is the Lie algebra denoted by $\f{n}_7$ in \cite{Andrada1}. The skew-symmetric bilinear map $\omega:\fg\times\fg\to {\mathbb R}$ defined by
$$\omega(E_1, E_3)=4, \ \omega(E_1, E_4)=-2, \ \omega(E_2,E_5)=-2, \ \omega(E_4,E_6)=4$$ and $\omega (E_j,E_k)=0$ for the other cases in which $j<k$, turns to be a symplectic form on $\fg$ and, further, the map $J\in\f{gl}(\fg)$ such that
$J(E_j)=E_{3+j}, \ J(E_{3+j})=-E_j$, for $j=1,2,3$,
is an abelian complex structure compatible with $\omega$. Therefore, $(\fg,\omega,J)$ is pseudo-K\"ahler. 

However, the left-symmetric product defined by $\omega$ is not associative. Actually, one can easily verify that, for instance, $(E_1\cdot E_4)\cdot E_1\ne E_1\cdot (E_4\cdot E_1)$. 

Note that this algebra provides an example of pseudo-K\"ahler Lie algebra with abelian complex structure which is not of the types of the two previous examples.
\end{enumerate}
}\end{examples}

A complex structure $J$ on a real Lie algebra $\fg$ provides a decomposition of the complexification $\fg^{\mathbb C}$ as the direct sum of 
the vector spaces $\fg^{1,0}=\nuc (J^{\mathbb C}-i\id),\fg^{0,1}=\nuc (J^{\mathbb C}+i\id)$, where $J^{\mathbb C}$ is the ${\mathbb C}$-linear map on  $\fg^{\mathbb C}$ defined by $J^{\mathbb C}(x_1+ix_2)=Jx_1+iJx_2$ for $x_1,x_2\in\fg$. One can easily verify that $\fg^{1,0}$ and $\fg^{0,1}$ are complex subalgebras of 
$\fg^{\mathbb C}$ and, obviously, they have the same dimension. Hence, the complex Lie algebra  $\fg^{\mathbb C}$ is naturally endowed with a paracomplex structure. Further, if the complex structure $J$ is abelian, then the subalgebras $\fg^{1,0}$ and $\fg^{0,1}$ are abelian and, thus, the paracomplex structure on $\fg^{\mathbb C}$ is also abelian \cite{BB-PK}. Further, if the real Lie algebra is pseudo-K\"ahler, then its complexification results to be para-K\"ahler. The following result will be used in next section. We omit its proof since it is nothing but a straghtforward calculation.

\smallskip

\begin{lemma} \label{prod_complexif} If  $(\fg, \omega, J)$ be a pseudo-K\"ahler Lie algebra with abelian complex structure then its complexification $\fg^{\mathbb C}$ is naturally endowed with the abelian para-K\"ahler structure defined by  the complex symplectic form $\omega^{\mathbb C}(x+iy,x'+iy')=\omega(x,x')-\omega(y,y')+i\omega(x,y')+i\omega(y,x')$, where $x,x',y,y'\in\fg$ and the abelian paracomplex structure $K(x+iy)=-Jy+iJx$. Further, the left-symmetric product defined by $\omega^{\mathbb C}$ on $\fg_{\mathbb C}$ is obtained from the left-symmetric product on $\fg$ as follows:
$$(x+iy)\cdot (x'+iy')=x\cdot x' - y\cdot y'+ ix\cdot y'+iy\cdot x' \, ,$$
for $x,y,x',y'\in\fg$\end{lemma}

\medskip

In order to fix notations, we will recall some well-known definitions on linear algebra. In the sequel, for a complex number $\alpha\in {\mathbb C}$, we will denote by $\mbox{Re}(\alpha)$ and $\mbox{Im}(\alpha)$ respectively its real and imaginary part and $\overline\alpha$ will mean its complex conjugate. 

\medskip

 Let $V$ be a complex vector space. A {\it semilinear map} $\tau:V\to V$ is a ${\mathbb R}$-linear map such that $\tau (\alpha v)=\overline{\alpha}\tau (v)$ for all $\alpha\in {\mathbb C}$ and $v\in V$. A {\it sesquilinear form} is a ${\mathbb R}$-bilinear map $H:V\times V\to {\mathbb C}$ which is ${\mathbb C}$-linear in the left component and semilinear in the right one. Further, a sesquilnear form $H$ is said to be {\it hermitian} if $H(v_1,v_2)=\overline{H(v_2,v_1)}$ holds for all $v_1,v_2\in V$. Note that, even though we will only consider non-degenerate hermitian forms, we shall not impose $H$ to be definite. 

Let $(V,H)$ be a complex vector space with a non-degenerate hermitian form. If $W$ is a complex subspace of $V$, then we will denote by $W^\bot$ its $H$-orthogonal  subspace, this is to say
$$W^\bot=\{\ x\in V\, :\, H(x,y)=0 \ \mbox{ for all }\ y\in W\ \}.$$
We remind that if $F:V\to V$ is a ${\mathbb C}$-linear map, its $H$-adjoint map $F^*:V\to V$ is uniquely defined by
$H(F^*x,y)=H(x,Fy)$ for $x,y\in V$.

\medskip

 \begin{rema} {\em All through the paper we will consider traces of ${\mathbb R}$-linear maps. Even for a ${\mathbb C}$-linear map, the trace considered will be the trace of the corresponding ${\mathbb R}$-linear map on the underlying real vector space. 

A simple calculation shows  that for a semilinear map $\tau$ on a complex vector space one always has $\mbox{trace}(\tau)=0$.}
\end{rema}

\section{Compatible complex associative commutative algebras}

As we will see below, every pseudo-K\"ahler Lie algebra with abelian complex structure is obtained from a pair $(U,H)$  where $U$ is a complex associative commutative algebra and $H$ a non-degenerate hermitian form on $U$ which verifies some compatibility conditions with respect to the associative product on $U$. We start this section with the construction of the associated pseudo-K\"ahler Lie algebras for a given pair $(U,H)$.

\begin{defi} {\em We will say that a complex associative commutative algebra $U$ endowed with a non-degenerate hermitian sesquilinear form $H$ is {\it  compatible with an abelian pseudo-K\"ahler structure} or, shortly, {\it APK-compatible} if and only if for every $x,y,z\in U$ it holds that
\begin{eqnarray}\label{compat} xR_{z}^{*}y=yR_{z}^{*}x\ ,\end{eqnarray}
where $R_z$ stands for the multiplication by $z$ in $U$.
}\end{defi}

\begin{prop} Let $U$ be a complex associative commutative algebra and $H$ a non-degenerate hermitian form such that the pair $(U, H)$ is APK-compatible.

Let  $\fg_U$ denote the underlying real vector space of $U$ and define $[\cdot,\cdot]:\fg_U\times\fg_U\to \fg_U$, $J:\fg_U\to\fg_U$ and $\omega :\fg_U\times\fg_U\to {\mathbb R}$ as follows:
$$
 [x,y]=R_{y}^*x-R_{x}^*y\ ,\quad
 Jx=i\, x\ ,\quad
\omega (x,y)=\mbox{\rm Im}(H(x,y))\ ,
$$
for all $x,y\in U$.

The pair $(\fg_U,[\cdot,\cdot])$ is a real Lie algebra and $(\omega,J)$ provides an abelian pseudo-K\"ahler structure on $\fg_U.$
\end{prop}
\dem 
First note that $R_{\alpha x}^*=\overline{\alpha}R_{x}^*$  and, obviously, $R_{x+y}^*=R_{x}^*+R_{y}^*$ for all $x,y\in U$ and $\alpha\in {\mathbb C}$. This shows that the bracket is ${\mathbb R}$-bilinear. Furhter, Jacobi identity follows immediately from the equalities
$$R_{[x,y]}^*=R_{y}R_{x}^*-R_{x}R_{y}^*\ ,\qquad 
R_{x}^*[y,z]=R_{zx}^*y-R_{yx}^*z\ ,$$
which can be easily derived from the condition given in equation (\ref{compat}) and the definition of the bracket.

 Since $H$ is non-degenerate so must be $\omega$ because for every $x,y\in U$ one has
$$H(x,y)=\omega (ix, y)+i\,\omega(x,y).$$
Further, $\omega$ is obviously skew-symmetric, because $H$ is hermitian, and for all $x,y,z\in U$ we get
\begin{eqnarray*}
& & H ([x,y],z)+H ([y,z],x)+H ([z,x],y)=\\ & &
H (R_{y}^*x-R_{x}^*y,z)+H (R_{z}^*y-R_{y}^*z,x)+H(R_{x}^*z-R_{z}^*x,y)=\\
& & H(x,yz)-H(y,xz)+H(y,zx)-H(z,yx)+H(z,xy)-H(x,zy)=0,
\end{eqnarray*}
and, therefore, its imaginary part also vanishes, which proves that $\omega$ is a symplectic form.

Finally, it is clear that $J$ turns to be an abelian complex structure on $\fg$ and that $\omega (Jx,Jy)=\omega (x,y)$. Thus $(\fg,\omega,J)$ is pseudo-K\"ahler.\qed

\begin{rema}\label{rema_apk} {\em The following facts concerning the pseudo-K\"ahler Lie algebra contructed above should be noticed:
\begin{enumerate}
\item[(1)] The pseudo-K\"ahler metric is given by
$g(x,y)=\omega (Jx,y)=\mbox{Re}(H(x,y))$
for all $x,y\in \fg_U$.
\item[(2)] The left-symmetric product defined by $\omega$ on $\fg_U$ is as follows:
$$x\cdot y= (R_y+R^*_y)x\, ,\quad x,y\in \fg_U.$$
\item[(3)]
 It is straightforward to show that when $U=U_1\oplus U_2$ is a direct sum of $H$-orthogonal ideals then $\fg_U$ is also the $\omega$-orthogonal direct sum of the ideals $\fg_{U_1}$ and $\fg_{U_2}$.\end{enumerate}
}
\end{rema}

The following result shows that every abelian pseudo-K\"ahler Lie algebra of type $\f{aff}({\mathcal A})$ can be viewed as the Lie algebra associated to a compatible pair $(U,H)$ and provides a necessary and sufficient condition to recognize such algebras among those constructed by APK-compatible pairs.

\begin{prop} \label{caracaff}
If $({\mathcal A},B)$ is a symmetric associative commutative algebra then $\f{aff}({\mathcal A})$ with the standard pseudo-K\"ahler structure is holomorphically symplectomorphic to the abelian pseudo-K\"ahler Lie algebra constructed  by the APK-compatible pair $(U,H)$ defined by the complexification $U={\mathcal A}^{\mathbb C}$ and the hermitian form given by $$H(a+ib,a'+ib')=-4B(a,a')-4B(b,b')-i4(B(b,a')-B(a,b'))$$ for all $a,b,a',b'\in {\mathcal A}$.

Conversely, if for an abelian pseudo-K\"ahler Lie algebra $\fg_U$ defined by an APK-compatible pair $(U,H)$ there exists a semilinear involution $\sigma: U\to U$ such that
$$\sigma (xy)=\sigma (x)\sigma (y)\,  ,\quad H(\sigma (x),\sigma (y))=\overline{H(x,y)}\, ,\quad H(xy,z)=H(x,\sigma(y)z)$$
hold for all $x,y,z\in U$, then $\fg_U$ is holomorphically symplectomorphic to  
$\f{aff}({\mathcal A})$ with its standard pseudo-K\"ahler structure, where ${\mathcal A}$ is the real form of $U$ defined by $\sigma$ endowed with a bilinear form $B$ obtained, up to constant, by the restriction of $H$ to ${\mathcal A}$.
\end{prop}
\dem Suppose first that $({\mathcal A},B)$ is a symmetric associative commutative algebra and define $H$ on $U={\mathcal A}^{\mathbb C}$ as in the statement. Notice that $H$ is non-degenerate if and only if $B$ is so. A simple calculation shows that $R^*_{a+ib}=R_{a-ib}$ for all $a,b\in {\mathcal A}$. Thus, the condition (\ref{compat}) follows immediately from the associativity and commutativity of $U$. If $\psi :\f{aff}({\mathcal A})\to \fg_U$ is defined by $\psi (a,b)=-(a+ib)/2$, one easily sees that $\psi$ commutes with the complex structures and that 
$\mbox{\rm Im}H(\psi (a,b),\psi (a',b'))=B(a,b')-B(b,a')$ for all $a,b,a',b'\in {\mathcal A}$. Moreover, $\psi$ is an isomorphism of real Lie algebras since it is clearly ${\mathbb R}$-linear and for $a,b,a',b'\in {\mathcal A}$ we have
\begin{eqnarray*}[\psi (a,b),\psi (a',b')]& =&\frac{1}{4}[a+ib,a'+ib']=\frac{1}{4}R^*_{a'+ib'}{(a+ib)}-\frac{1}{4}R^*_{a+ib}{(a'+ib')}\\&=&
\frac{1}{4}R_{a'-ib'}{(a+ib)}-\frac{1}{4}R_{a-ib}{(a'+ib')}=
\frac{i}{2}(a'b-b'a)=\psi (0,ab'-ba')\\&=&\psi ([(a,b),(a',b')].\end{eqnarray*}

Conversely, if $(U,H)$ is an APK-compatible algebra admitting a semilinear involution $\sigma$ verifying the properties given above, then $U$ admits the associative commutative real form
${\mathcal A}=\{x+\sigma (x)\, ;\, x\in U\}$. For $x,y,z\in U$ we have
\begin{eqnarray*}
& & H(x+\sigma (x),y+\sigma (y))=2\ \mbox{\rm Re}(H(x,y)+H(\sigma (x),y))\\
& & H((x+\sigma (x))(y+\sigma (y)),z+\sigma (z))=H(x+\sigma (x),(y+\sigma (y))(z+\sigma (z))).
\end{eqnarray*}
From the first identity we get that the form $B$ defined on ${\mathcal A}$ by $B(a,b)=-\frac{1}{4}H(a,b)$ is  a real-valued  symmetric bilinear form on ${\mathcal A}$ and, from the second, that $B(ab,c)=B(a,bc)$ for $a,b,c\in {\mathcal A}$. Therefore, $({\mathcal A},B)$ is a symmetric associative commutative algebra and, according to the first part of the proposition, $\f{aff}({\mathcal A})$ is holomorphically symplectomorphic to the Lie algebra constructed with the APK-compatible pair $({\mathcal A}^{\mathbb C},\tilde{H})$ where 
$\tilde{H}(a+ib,a'+ib')=-4B(a,a')-4B(b,b')-i4(B(b,a')-B(a,b'))$ for all $a,b,a',b'\in {\mathcal A}$. Since $U={\mathcal A}^{\mathbb C}$, we only need to prove that $\tilde{H}=H$. But this follows at once since both
$\tilde{H}$ and $H$ are hermitian and coincide on the real form ${\mathcal A}$ of $U$.\qed

\begin{lemma} \label{difference} If $(U,H)$ is an APK-compatible pair and  $\mbox{\rm ann}(U)$ denotes the annihilator of $U$, then 
$$(R_{x}^*R_z-R_zR_{x}^*)y\in \mbox{\rm ann}(U)$$
for all $x,y,z\in U$.\end{lemma}
\dem 
Let us consider $u,x,y,z\in U$. Since the product in $U$ is associative and commutative, we have
$$ u\left(R_{x}^*(zy)-zR_{x}^*y\right)=uR_{x}^*(zy)-uzR_{x}^*y=zyR_{x}^*u-z\left(uR_{x}^*y\right)=
 yzR_{x}^*u-zyR_{x}^*u=0\ ,$$
which proves that $R_{x}^*(zy)-zR_{x}^*y$ is in the annihilator.\qed

\begin{lemma} Let $(U,H)$ be an APK-compatible pair and let $\mbox{\rm ann}(U)$ denote the annihilator of $U$ and $R_U^*U$ the linear ${\mathbb C}$-span of all the elements $R_{x}^*y$ with $x,y\in U$. The following holds:
\begin{enumerate}
\item[(a)] The vector space $R_U^*U$ is the $H$-orthogonal subspace of  $\mbox{\rm ann}(U)$ and it is an ideal of the associative commutative algebra $U$.
\item[(b)] If $\mbox{\rm ann}(U)=\{0\}$ then $U=U^2=R_U^*U$.
\end{enumerate}
\end{lemma}
\dem It is clear that $H(R_x^*y,x_0)=0$ for all $x,y\in U$ if and only if $H(y,xx_0)=0$, this is to say, $x_0\in  \mbox{\rm ann}(U)$. In order to prove that $R_U^*U$ is an ideal, take $x,y,z\in U$ and $x_0\in  \mbox{\rm ann}(U)$. We then have
$$H(xR_{y}^*z,x_0)=H(z,yR_{x}^*x_0)=H(z,x_0R_{x}^*y)=0,$$
which proves that $xR_{y}^*z\in \mbox{\rm ann}(U)^\bot=R_U^*U$. This completes the proof of (a).

Now, if  $\mbox{\rm ann}(U)=\{0\}$ then $R_U^*U=\mbox{\rm ann}(U)^\bot=U$ and hence it only remains to see that it coincides with $U^2$.
Let us  take an element $x\in (U^2)^\bot$. For all $y,z,t\in U$ we then have
$$H(t,xR_{y}^*z)=H(t,zR_{y}^*x)=H(yR_z^*t,x)=0.$$
This clearly implies that $xR_{y}^*z=0$ for all $y,z\in U$ but, since $U= R_U^*U$, one obtains $x\in \mbox{\rm ann}(U)=\{0\}$. Hence, the $H$-orthogonal subspace of $U^2$ is null, which shows that $U=U^2$.\qed

\begin{prop} \label{annzero} Let $(U,H)$ be an APK-compatible pair and $(\fg_U,\omega,J)$ the corresponding pseudo-K\"ahler Lie algebra with abelian complex structure.

If $\mbox{\rm ann}(U)=\{0\}$, then  $(\fg_U,\omega,J)$ is holomorphically symplectomorphic to $\f{aff}({\mathcal A})$ for some real symmetric associative commutative algebra $({\mathcal A},B)$  endowed with the standard pseudo-K\"ahler structure.
\end{prop}
\dem If $\mbox{\rm ann}(U)=\{0\}$ we have, according to the lemma above, that $U=R^*_UU$. We can then define a map $\sigma :U\to U$ by $\sigma \left(R^*_xy\right)=R^*_yx$. Note that $\sigma$ is well defined since if $R^*_xy=R^*_zt$ then for all $u,v\in U$ we get
$$H(R^*_yx-R^*_tz,R^*_uv)=H(R^*_yx,R^*_uv)-H(R^*_tz,R^*_uv)=H(R^*_vu,R^*_xy)-H(R^*_vu,R^*_zt)=0.$$
It is obvious that $\sigma^2=\mbox{id}_U$ and that $$\sigma (\alpha R^*_xy)=\sigma ( R^*_x(\alpha y))=R^*_{\alpha y}x=\overline{\alpha}R^*_yx=\overline{\alpha}\sigma ( R^*_xy).$$
Further, using Lemma \ref{difference}, we get that $R_{x}^*(yz)=z\left(R_{x}^*y\right)$ holds for all $x,y,z\in U$ and, therefore, 
$$\sigma \left((R^*_xy)(R^*_zt)\right)=\sigma \left(tR^*_z(R^*_xy)\right)=\sigma \left(R^*_{zx}(yt)\right)=R^*_{yt}(xz)=\left(R^*_{y}x\right)\left(R^*_{t}z\right)=
\sigma \left(R^*_xy\right)
\sigma \left(R^*_zt\right).$$
Bearing in mind Proposition \ref{caracaff} it only remains to show that the identities $H(\sigma (x),\sigma (y))=\overline{H(x,y)}$ and $H(xy,z)=H(x,\sigma(y)z)$ are verified for all $x,y,z\in U$. But this follows at once because from equation (\ref{compat}) we have
\begin{eqnarray*} & & H(R^*_uv,R^*_xy)=H(R^*_yx,R^*_vu)=\overline{H(R^*_vu,R^*_yx)}\\
& & H(xR^*_uv,z)=H(R^*_uv,R^*_xz)=H(R^*_ux,R^*_vz)=H(x,uR^*_vz)=H(x,zR^*_vu),\end{eqnarray*}
which yields the desired identity.\qed

\bigskip

 In the cases where $\mbox{ann}(U)\cap (\mbox{ann}(U))^\bot =\{0\}$ something similar occurs, as we prove in the following proposition which also clarifies the structure of the corresponding real symmetric associative commutative algebra. 

\
\begin{prop} \label{aff_ext} Let $(U,H)$ be an APK-compatible pair. If $\mbox{\rm ann}(U)\cap (\mbox{\rm ann}(U))^\bot =\{0\}$ then $\fg_U$ is holomorphically symplectomorphic to an algebra $\f{aff}({\mathcal A})$ endowed with the standard pseudo-K\"ahler structure for a certain symmetric form $B$.

Moreover, the algebra ${\mathcal A}$ decomposes as a $B$-orthogonal sum of ideals 
$${\mathcal A}=\mbox{\rm ann}({\mathcal A})\oplus {\mathcal A}_1\oplus\cdots \oplus {\mathcal A}_p$$
where $p$ is the number of simple ideals in the semisimple part of ${\mathcal A}$ and each ${\mathcal A}_i$ is a unital algebra.
\end{prop}
\dem If $\mbox{\rm ann}(U)\cap (\mbox{ann}(U))^\bot =\{0\}$ then $U=\mbox{ann}(U)\oplus (\mbox{ann}(U))^\bot$, which is an orthogonal sum of ideals. The pseudo-K\"ahler Lie algebra constructed with the algebra
$\mbox{ann}(U)$ is clearly abelian and, therefore, holomorphically symplectomorphic to $\f{aff}({\mathcal A}_0)$ for a real algebra with ${\mathcal A}_0^2=\{0\}$.

The ideal $\tilde{U}=(\mbox{ann}(U))^\bot$ obviously verifies $\mbox{ann} (\tilde{U})=\{0\}$ and according to Proposition \ref{annzero} $\fg_{\tilde{U}}$ is holomorphically symplectomorphic to $\f{aff}(\tilde{\mathcal A})$ for some real symmetric associative commutative algebra $(\tilde{\mathcal A},\tilde{B})$. 
 From $\mbox{ann} (\tilde{U})=\{0\}$ we immediately get that $\mbox{ann}(\tilde{\mathcal A})$ is also null. This implies that $\tilde{\mathcal A}$ cannot be nilpotent and hence it has a non trivial semisimple part. Since $\tilde{\mathcal A}$ is commutative, its semisimple part is the direct sum of several copies of ${\mathbb R}$ and ${\mathbb C}$ with their usual field structure. Let $p$ be the number of those simple ideals and let us choose on each of them an idempotent element $e_i$. Consider ${\mathcal A}_i=R_{e_i}\tilde{\mathcal A}$. It is clear that ${\mathcal A}_i$ is an ideal of $\tilde{\mathcal A}$ for all $i\leq p$ with unity $e_i$, that ${\mathcal A}_i{\mathcal A}_j=\{0\}$ and that they are mutually orthogonal since when $i\neq j$ we have $B(e_ia,e_jb)=B(a,e_ie_jb)=0$ for all $a,b\in\tilde{\mathcal A}.$ To see that $\tilde{\mathcal A}={\mathcal A}_1\oplus\cdots \oplus {\mathcal A}_p$, take the Peirce decomposition \cite{albert} $\tilde{\mathcal A}={e_1}\tilde{\mathcal A}\oplus {\mathcal M}_1$ where ${\mathcal M}_1=\{ a\in \tilde{\mathcal A}\ :\ e_1a=0\}$. Obviously, ${\mathcal A}_2+\cdots + {\mathcal A}_p\subset {\mathcal M}_1$ and ${\mathcal M}_1$ has trivial annhilator. We can proceed with ${\mathcal M}_1$ in the same way, taking its Pierce decomposition relative to the idempotent $e_2$. Repeating sucessively the same argument, we obtain that 
$\tilde{\mathcal A}={\mathcal A}_1\oplus\cdots \oplus {\mathcal A}_p\oplus {\mathcal M}_p$ where ${\mathcal M}_p=\{ a\in \tilde{\mathcal A}\ :\ e_ia=0, \mbox{ for all }i\le p\}$ is a nilpotent ideal of $\tilde{\mathcal A}$. If we suppose that ${\mathcal M}_p$ is not null, then we may find an element $a\in\mbox{ann}( {\mathcal M}_p),a\ne 0$. But this would imply $a\in \mbox{ann}(\tilde{\mathcal A})$, a contradiction. Hence ${\mathcal M}_p=\{0\}$ and $\tilde{\mathcal A}={\mathcal A}_1\oplus\cdots \oplus {\mathcal A}_p$.

Finally, if we take the orthogonal sum ${\mathcal A}={\mathcal A}_0\oplus\tilde{\mathcal A}$  one easily proves that $\mbox{ann}({\mathcal A})={\mathcal A}_0$ and that 
$\f{aff}({\mathcal A})=\f{aff}({\mathcal A_0})\oplus\f{aff}(\tilde{\mathcal A})$, and the result follows.\qed

\begin{rema}{\em  A well-known result on unitary symmetric algebras implies that and the restriction $B_i$  of $B$ to ${\mathcal A_i}\times {\mathcal A_i}$ is given by $B_i(a,b)=f_i(ab)$ for all $a,b\in {\mathcal A}_i$ where $f_i:{\mathcal A}_i\to {\mathbb R}$ is a linear form such that $f_i\circ R_a\neq 0$ for all $a\in {\mathcal A}_i$.}
\end{rema}

\bigskip

The main result of this section shows that every abelian pseudo-K\"ahler Lie algebra is associated to certain APK-compatible complex associative commutaive algebra. 

\medskip

\begin{theor} Let $(\fg, \omega, J)$ be a pseudo-K\"ahler Lie algebra with abelian complex structure. There exists an APK-compatible pair $(U,H)$ such that $\fg$ and $\fg_U$ are holomorphically symplectomorphic.
\end{theor}
\dem If $\fg$ is pseudo-K\"ahler, then  $(\fg^{\mathbb C},\omega^{\mathbb C},K)$ is a para-K\"ahler complex Lie algebra for $\omega^{\mathbb C}$ and $K$ defined as in Lemma \ref{prod_complexif}. Note that $\fg^{1,0}=\mbox{ker}(J^{\mathbb C}-i\mbox{id})=\mbox{ker}(K+\mbox{id})$ and,  according to \cite[Lemma 4.1]{BB-PK}, it is a complex associative commutative algebra. Let us then consider $$U=\fg^{1,0}=\{x-iJx\, :\, x\in\fg\}.$$ Recall that the associative product on $U$ is nothing but the restriction of the left-symmetric product defined on $\fg^{\mathbb C}$ by $\omega^{\mathbb C}$, this is to say, for $x,y,z,t\in\fg$ we have
$$\omega^{\mathbb C}((x-iJx)(y-iJy),z+it)=-\omega^{\mathbb C}(y-iJy,[x-iJx,z+it]).$$
Now we can define the non-degenerate hermitian form $H:U\times U\to {\mathbb C}$  as follows:
$$H(x-iJx,x'-iJx')={2i}\omega^{\mathbb C}(x-iJx,x'+iJx')=-4\omega(x,Jx')+i4\omega(x,x'),$$ for $x,x'\in\fg$. We will now show that
$u'R^*_{u}u''=u''R^*_{u}u'$ for all $u,u',u''\in U$. Let us consider
$$\overline{U}=\fg^{0,1}=\{x+iJx\, :\, x\in\fg\}$$
and, for $x,y\in \fg$, let us put $\overline{x+iy}=x-iy$. Since $\omega ^{\mathbb C}(\overline{\xi_1},\overline{\xi_2})=\overline{
\omega ^{\mathbb C}({\xi_1},{\xi_2})}$ and $[\overline{\xi_1},\overline{\xi_2}]=\overline{
[{\xi_1},{\xi_2}]}$ for all $\xi_1,\xi_2\in \fg^{\mathbb C}$, for the left-symmetric product  defined by $\omega^{\mathbb C}$ one also has $\overline{\xi_1}\cdot\overline{\xi_2}=\overline{\xi_1\cdot\xi_2}$ because 
$$\omega^{\mathbb C}(\overline{\xi_1\cdot\xi_2},\xi)=\overline{\omega^{\mathbb C}(\xi_1\cdot\xi_2,\overline{\xi})}=
-\overline{\omega^{\mathbb C}(\xi_2,[\xi_1,\overline{\xi}])}=-\omega^{\mathbb C}(\overline{\xi_2},[\overline{\xi_1},{\xi}])=
\omega^{\mathbb C}(\overline{\xi_1}\cdot\overline{\xi_2},{\xi})$$
holds for all $\xi\in\fg^{\mathbb C}$. Further, one easily verifies that  $\fg^{\mathbb C}=U\oplus\overline{U}$ verifies
$$\omega^{\mathbb C}(U,U)=\omega^{\mathbb C}(\overline{U},\overline{U})=\{0\}\, ,\quad u\cdot \overline{u'}\in U\, ,\quad \overline{u}\cdot u\in \overline{U}\, ,  $$
for all $u,u'\in U$ (which can also be deduced from the results on para-K\"ahler algebras given in \cite{BB-PK}). Thus, for all $u,u',u''\in U$, we have
$$H(R^*_uu',u'')=H(u',uu'')=2i\omega^{\mathbb C}(u',\overline{u}\overline{u''})=
-2i\omega^{\mathbb C}([\overline{u},u'],\overline{u''})=
2i\omega^{\mathbb C}(u'\cdot\overline{u},\overline{u''})=H(u'\cdot\overline{u},{u''}).$$
Therefore we get that $R_u^*u'=u'\cdot \overline{u}$ for all $u,u'\in U$. Since $U$ is commutative and $(\fg^{\mathbb C},\cdot)$ is left-symmetric we get
$$0=[u',u'']\cdot \overline{u}=u' (u''\cdot\overline{u})-u'' (u'\cdot\overline{u})=u'R_u^*u''-u''R_u^*u'\, \quad u,u',u''\in U.$$
We have, thus, proved that $(U,H)$ is an APK-compatible pair. 

To see that $\fg$ and $\fg_U$ are holomorphically symplectomorphic, we will consider the map $\psi:\fg\to\fg_U$ given by $\psi (x)=\frac{1}{2}(x-iJx)$ for all $x\in\fg$. It is clear that $\psi (Jx)=i\psi(x) $, that $\omega (x,y)=\mbox{Im}(H(\psi(x),\psi(y))$
and that $\psi$ is ${\mathbb R}$-linear. Let us prove that $\psi[x,y]=[\psi(x),\psi(y)]_{\fg_U}$. But, if we take $x,y\in\fg$, then
\begin{eqnarray*}
 4[\psi(x),\psi(y)]_{\fg_U}&=&4R^*_{\psi (y)}\psi (x)-4R^*_{\psi (x)}\psi (y)=4\psi(x)\cdot \overline{\psi(y)}-4\psi(y)\cdot \overline{\psi(x)}\\&=&
(x-iJx)\cdot (y+iJy)-(y-iJy)\cdot (x+iJx)\\
& =& x\cdot y+Jx\cdot Jy-iJx\cdot y+ix\cdot Jy-y\cdot x-Jy\cdot Jx+iJy\cdot x-iy\cdot Jx\\&=&[x,y]+[Jx,Jy]-iJ[x,y]-iJ[Jx,Jy]=2([x,y]-iJ[x,y])=4\psi[x,y].
\end{eqnarray*}
Since $\psi$ is obviously invertible, it is an isomorphism of real Lie algebras.\qed

\bigskip

We shall use that every Lie algebra with abelian pseudo-K\"ahler structure is of the form $\fg_U$ for some APK-compatible pair $(U,H)$ to see when such a Lie algebra is unimodular. We first prove the following lemma:

\begin{lemma} Let $(U,H)$ be an APK-compatible pair and consider $x\in U$. For each $k\in {\mathbb N}$, $k\geq 2$ then the power $(R_x+R^*_x)^k$ is a linear combination of 
$R_x^k$, $(R^*_x)^k$, $R_x^r(R^*_x)^{k-r}$ and $(R^*_x)^{k-r}R_x^r$, $1\le r \le k-1$.
\end{lemma}
\dem We will proceed by induction on $k$. The result is clear for $k=2$ because
$$(R_x+R^*_x)^2=R_x^2+R_xR^*_x+R^*_xR_x+(R^*_x)^2.$$

From Lemma \ref{difference} we obtain
 that $R_xR^*_yR_x=R_x^2R^*_y$ and thus $R_x^*R_yR^*_x=(R_xR^*_yR_x)^*=R_y(R^*_x)^2$. Therefore, when $s,r\ge 1$ we have
$$R^*_xR_x^r(R^*_x)^s=R_x^r(R^*_x)^{s+1}\, ,\quad R_x(R^*_x)^sR_x^r=R_x^{r+1}(R^*_x)^s\, .$$

Now, we have that $(R_x+R^*_x)R^{k-1}_x=R^k_x+R^*_xR^{k-1}$, $(R_x+R^*_x)(R^*_x)^{k-1}=R_x(R^*_x)^{k-1}+(R^*_x)^{k}$ and for $s,r\geq 1$ we get
\begin{eqnarray*}& & (R_x+R^*_x)R_x^r(R^*_x)^s=R_x^{r+1}(R^*_x)^s+R^*_xR_x^r(R^*_x)^s=R_x^{r+1}(R^*_x)^s+R_x^r(R^*_x)^{s+1}\, ,\\
& & (R_x+R^*_x)(R^*_x)^sR_x^r=R_x(R^*_x)^sR_x^r+(R^*_x)^{s+1}R_x^r=R_x^{r+1}(R^*_x)^s+(R^*_x)^{s+1}R_x^r\, .
\end{eqnarray*}
and this shows that if $(R_x+R^*_x)^{k-1}$ verifies the condition, so does $(R_x+R^*_x)^k$.\qed

\begin{prop} Let $\fg$ be a pseudo-K\"ahler Lie algebra with abelian complex structure and $(U,H)$ an APK-compatible pair such that $\fg=\fg_U$. The following conditions are equivalent:
\begin{enumerate}
\item[(i)] $U$ is nilpotent,
\item[(ii)] $\fg$ is nilpotent,
\item[(iii)] $\fg$ is unimodular.
\end{enumerate}
\end{prop}
To see that $\fg$ is nilpotent whenever $U$ is so, let us consider $x,y\in U$. Recall that $R^*_{[x,y]}=R_yR^*_x-R_xR^*_y$ and hence
we have
$$\ad (x)[x,y]= yR^*_xx-xR^*_yx-R^*_x[x,y]=xR^*_xy-xR^*_yx-R^*_x[x,y]=-(R_x+R^*_x)[x,y]\,,$$
and, therefore, $\ad(x)^{k+1}(y)=(-1)^k(R_x+R^*_x)^k[x,y]$ for all $k\geq 1$. If we have that $U^n=\{0\}$ then the Lemma above implies that $(R_x+R^*_x)^{2n}=0$ and, thus, $\ad (x)^{2n+1}=0$, which proves the nilpotency of  $\fg$. 

It is obvious that if $\fg$ is nilpotent then it is unimodular and, thus, it only remains to show that for a unimodular $\fg$ we get $U$ nilpotent. Note that the map $G_x=\mbox{ad}(x)+R^*_x$ is semilinear and, therefore, traceless. This shows that $\mbox{trace}(\ad (x))=-\mbox{trace}(R^*_x)=-\mbox{trace}(R_x)$. Consequently, if $\fg$ is unimodular, then we have $\mbox{trace}(R_x)=0$ for all $x\in U$. But this implies that $U$ is nilpotent since, otherwise, there must exist an idempotent element $e\in U$ and one easily proves using the Peirce decomposition of $U$ that $\mbox{trace}(R^*_e)=d$ where $d$ is the real dimension of the ideal $eU$.\qed

\section{Curvature of the pseudo-K\"ahler metric}

We will now calculate the curvatures of the pseudo-Riemannian metric of a Lie algebra with abelian pseudo-K\"ahler  structure. We shall give two different approaches: a first one using the left-symmetric product defined by the symplectic form, and a second one based on the description in terms of APK-compatible pairs. 

Let us consider a Lie algebra with abelian pseudo-K\"ahler structure $(\fg,\omega,J)$ and denote by $(\fg,\cdot)$ the left-symmetric algebra structure defined by $\omega$.

\begin{prop} \label{levi-curv} Let $(\fg,J,\omega)$ be a pseudo-K\"ahler Lie algebra with abelian complex structure and $g$ the pseudo-K\"ahler metric. The Levi-Civita connection and the Riemannian curvature tensor of $g$ are respectively given for $x,y,z\in\fg$ by
$$\nabla_xy=-Jy\cdot Jx\, ,\qquad \mathbf{R}(x,y)z=-Jz\cdot J[x,y]+(z\cdot Jy)\cdot Jx-(z\cdot Jx)\cdot Jy\, .$$

If ${\mathcal K}$ denotes the Killing form of $\fg$, then the Ricci curvature tensor  $\mbox{\rm Ric}(x,y)=\mbox{\rm trace}\Big(\mathbf{R}(x,-)y\Big)$ is given by
$$\mbox{\rm Ric}(x,y)=\mbox{\rm trace}(\mbox{\rm ad}(Jy\cdot Jx))-{\mathcal K}(x,y)\, ,\quad x,y\in\fg.$$
\end{prop}
\dem
 For $x,y,z\in\fg$ the Levi-Civita connection is given by
\begin{eqnarray*}
 2g(\nabla_xy,z)&=&g ([x,y],z)-g ([y,z],x)+g ([z,x],y) =\omega (J[x,y],z)-\omega (J[y,z],x)+\omega (J[z,x],y) \\
&= & \omega (J[x,y],z)-\omega (z,y\cdot Jx)-\omega (z,x\cdot Jy)  = \omega (J[x,y]+y\cdot Jx+x\cdot Jy,z)\\&= &g ([x,y]-J(y\cdot Jx)-J(x\cdot Jy),z)
\end{eqnarray*}
and therefore, since $J$ commutes with right multiplications and $[x,y]=[Jx,Jy]$, we have
$$2\nabla_xy=[Jx,Jy]-Jy\cdot Jx-Jx\cdot Jy=-2Jy\cdot Jx,$$
for all $x,y\in\fg$.

Now, the curvature tensor is  given by
\begin{eqnarray*}
 \mathbf{R}(x,y)z&=&\nabla_{[x,y]}z-\nabla_x\nabla_yz+\nabla_y\nabla_xz =-Jz\cdot J[x,y]-J(Jz\cdot Jy)\cdot Jx+J(Jz\cdot Jx)\cdot Jy\\
& =& -Jz\cdot J[x,y]+(z\cdot Jy)\cdot Jx-(z\cdot Jx)\cdot Jy
\end{eqnarray*}
for all $x,y,z\in\fg$.

In order to compute the Ricci form, let  ${\mathcal L_x}$ denote the left multiplication by an element $x\in\fg$ with the left-symmetric product induced by $\omega$. Then we have $x'\cdot Jx=(\mbox{ad}(x)J+{\mathcal L}_{Jx})x'$ for all $x'\in\fg$ and it follows that
$$\mathbf{R}(x,z)y=-{\mathcal L}_{Jy}J\mbox{ad}(x)z+(\mbox{ad}(x)J+{\mathcal L}_{Jx}){\mathcal L}_yJz-{\mathcal L}_{y\cdot Jx}Jz.$$
Recalling that $J{\mathcal L}_{y}={\mathcal L}_{Jy}$, we get
\begin{eqnarray*} \mbox{Ric}(x,y)&=&-\mbox{trace}({\mathcal L}_{Jy}J\mbox{ad}(x))+\mbox{trace}(\mbox{ad}(x){\mathcal L}_{Jy}J)+\mbox{trace}({\mathcal L}_{Jx}{\mathcal L}_yJ)-\mbox{trace}({\mathcal L}_{y\cdot Jx}J)\\
&=& -\mbox{trace}({\mathcal L}_{x}{\mathcal L}_y)-\mbox{trace}({\mathcal L}_{Jy\cdot Jx}).\end{eqnarray*}
Now, the result follows immediately from the fact that for all $z\in\fg$, the map ${\mathcal L}_z$ is the adjoint map with respect to $\omega$ of $-\mbox{ad}(z)$ and that a linear map and its $\omega$-adjoint have the same trace.\qed

\begin{rema}{\em Observe that, actually, Lemma \ref{commutes} is equivalent to the parallelism of $J$ with respect to the Levi-Civita connection.}
\end{rema}

\begin{example}\label{curv_aff} {\em 
Let $ ({\mathcal A},B)$ be a real symmetric associative commutative algebra  and let us consider on $\f{aff}({\mathcal A})$ the standard pseudo-K\"ahler structure. We had already seen that the left-symmetric product defined by the symplectic form was given by
$(a,b)\cdot (a',b')=-(aa',ba')$ for all $a,b,a',b'\in {\mathcal A}$. Consequently, the levi-Civita connection and the Riemannian curvature tensor are
$$\nabla_{(a,b)}(a',b')=(b'b,-a'b)\, ,\quad \mathbf{R}((a,b),(a',b'))(a'',b'')=(ab'b''-ba'b'',ba'a''-ab'a'')\, ,\,\,\, a,b,a',b',a'',b''\in {\mathcal A}.$$
Note that this clearly implies that the metric is flat if and only if ${\mathcal A}^3=\{0\}$

Besides, we have that
$$\mathbf{R}((a,b),(a',0))(a'',b'')=(-R_{bb''}a',R_{ba''}a')\, ,\quad \mathbf{R}((a,b),(0,b'))(a'',b'')=(R_{ab''}b',-R_{aa''}b')\, ,$$
from where it is clear that
$$\mbox{Ric}((a,b),(a'',b''))=-\mbox{trace}(R_{bb''})-\mbox{trace}(R_{aa''}).$$
Then, the metric is Ricci flat if and only if $\mbox{trace}(R_{aa'})=0$ for all $a,a'\in {\mathcal A}$. But this only occurs if and only if ${\mathcal A}$ is nilpotent because, otherwise, there exists an idempotent  $e\in{\mathcal A}$ and one always has $\mbox{trace}(R_{ee})=\mbox{dim}(e{\mathcal A})\ge 1$.
}\end{example}

Flatness or Ricci flatness are easier to study if one describes the Riemannian and the Ricci curvature in terms of the products of an APK-compatible pair. We have the following:

\begin{prop} Let $\fg_U$ denote the pseudo-K\"ahler Lie algebra with abelian complex structure defined by an APK-compatible pair $(U,H)$. The Levi-Civita connection, the Riemannian curvature tensor and the Ricci curvature are given for $x,y\in U$ by
$$
\nabla_xy=(R_x-R_x^*)y\, ,\quad
\mathbf{R}(x,y)=3R_xR_y^*-3R_yR^*_x+R^*_xR_y-R^*_yR_x\, ,\quad
 \mbox{\rm Ric}(x,y)=-2\,\mbox{\rm trace}(R_yR^*_x)\, .
$$
\end{prop}
\dem According to (2) in Remark \ref{rema_apk}, and Proposition \ref{levi-curv}, we get that the Levi-Civita connection is
$$\nabla_xy=-Jy\cdot Jx=-R_{ix}iy-R^*_{ix}iy=R_xy-R^*_xy\, .$$
Thus, $\nabla_{[x,y]}=R_{[x,y]}-R^*_{[x,y]}=2R_xR_y^*-2R_yR^*_x$ and we then have
$$
\mathbf{R}(x,y)=2R_xR_y^*-2R_yR^*_x-(R_x-R_x^*)(R_y-R_y^*)+(R_y-R_y^*)(R_x-R_x^*)=3R_xR_y^*-3R_yR^*_x+R^*_xR_y-R^*_yR_x\, .$$

In order to compute the Ricci curvature, let us consider for each $x\in\fg$ the semilinear map $G_x:U\to U$ defined by $G_x(y)=R^*_yx$. It follows that
$$\mathbf{R}(x,z)y=(3R_xG_y-3R_yR^*_x+R^*_xR_y-G_{xy} )(z)\, .$$
Since $R_xG_y$ and $G_{xy}$ are semilinear, they are traceless and as a consequence we get
$$\mbox{Ric}(x,y)=-3\mbox{trace}(R_yR^*_x)+\mbox{trace}(R_xR^*_y)=-2\mbox{trace}(R_yR^*_x)\, ,$$
as claimed.\qed

\begin{coro} Let $(U,H)$ be an APK-compatible pair. The pseudo-K\"ahler metric on $\fg_U$ is flat if and only if $\ 3R_xR^*_x=R_x^*R_x$ for all $x\in U$.
\end{coro}
\dem If  $3R_xR^*_x-R_x^*R_x=0$ holds for all $x\in U$ then we have
\begin{eqnarray*}
& & 0=3R_{x+y}R^*_{x+y}-R^*_{x+y}R_{x+y}=3R_xR^*_y+3R_yR^*_x-R^*_xR_y-R^*_yR_x\, ,\\
& & 0=3R_{x+iy}R^*_{x+iy}-R^*_{x+iy}R_{x+iy}=-3iR_xR^*_y+3iR_yR^*_x-iR^*_xR_y+iR^*_yR_x\, ,\end{eqnarray*}
which imply that $3R_xR^*_y=R_y^*R_x$  for $x,y\in U$, and one obviously deduces that $\mathbf{R}(x,y)=0$. 

Conversely, if the metric is flat we have
$$0=R(x,ix)=-6iR_xR^*_x+2iR^*_xR_x\, $$
which gives the desired identity.\qed

\begin{prop} \label{flatU} Let $(U,H)$ be an APK-compatible pair and $(\fg_U,\omega, J)$ the associated pseudo-K\"ahler Lie algebra. Let us consider the following conditions:
\begin{enumerate}
\item[(c1)] The pseudo-K\"ahler metric is flat.
\item[(c2)] The left-symmetric product defined by $\omega$ is, actually, associative.
\item[(c3)] The derived ideal $[\fg,\fg]$ is contained in $\mbox{\rm ann}(U).$
\end{enumerate}
If one of the conditions above is fulfilled, then the other two are equivalent.
\end{prop}
\dem Recall that, as in the proof of the corollary above, the metric is flat if and only if $3R_xR^*_y=R_y^*R_x$ for all $x,y\in U$. On the other hand, the left-symmtric product is associative if and only if 
$$0=R_{x\cdot y}+R^*_{x\cdot y}-(R_y+R^*_y)(R_x+R^*_x)=R_xR^*_y-R^*_yR_x\, ,\quad x,y\in U.$$
Finally, $[\fg,\fg]\subset \mbox{\rm ann}(U)$ is equivalent to $R_xR^*_y=0$ for all $x,y\in U$ because we have
\begin{eqnarray*}
& & z[x,y]=zR^*_yx-zR^*_xy=(R_xR^*_y-R_yR^*_x)z
\\
& & z[ix,y]=izR^*_yx+izR^*_xy=i(R_xR^*_y+R_yR^*_x)z
\end{eqnarray*}
 for all $x,y,z\in U$. Now, it suffices to realize that the combination of two of these conditions imply $R_xR^*_y=R^*_xR_y=0$ and, then, the remaining condition is also verified.\qed

\begin{rema} {\em A left-symmetric algebra $(\f{A},\cdot)$ is called a {\it Novikov algebra} if  $(a\cdot b)\cdot c=(a\cdot c)\cdot b$ is also verified for all $a,b,c\in {\f{A}}$ \cite{bai},\cite{burde}. If $\fg$ is a pseudo-K\"ahler Lie algebra constructed with an APK-compatible pair $(U,H)$ and $(\fg,\cdot)$ denotes the left-symmetric algebra defined by the symplectic form then we have
$$(x\cdot y)\cdot z-(x\cdot z)\cdot y=(R_z+R^*_z)(R_y+R^*_y)x-(R_y+R^*_y)(R_z+R^*_z)x=(R_zR^*_y+R^*_zR_y-R_yR^*_z-R^*_yR_z)x$$
and, accordingly,
$$(x\cdot y)\cdot iz-(x\cdot iz)\cdot y=(iR_zR^*_y-iR^*_zR_y+iR_yR^*_z-iR^*_yR_z)x.$$
From these two equations we immediately get that $(\fg,\cdot)$ is Novikov if and only if $R_zR^*_y=R^*_yR_z$. This proves that, for our algebras, the Novikov condition is equivalent to the associativity of the left-symmetric product. 
}\end{rema}

\begin{coro} \label{flatU2} If a pseudo-K\"ahler algebra  with abelian complex structure $(\fg,\omega, J)$ is 2-step nilpotent, then the pseudo-K\"ahler metric is flat if and only if the left-symmetric product defined by $\omega$ is associative.
\end{coro} 
\dem Let $(U,H)$ be an  APK-compatible pair for $\fg$. A simple computation shows that for all $x,y,z\in U$ one has $[z,[x,y]]=-(R_z+R^*_z)[x,y]$ and, therefore, if $\fg$ is 2-nilpotent, then for $x,y,z\in U$ one has
$$0=[Jz,[x,y]]+J[z,[x,y]]=-(iR_z-iR^*_z)[x,y]-i(R_z+R^*_z)[x,y]=-2iz[x,y]\, ,$$
which shows that $[\fg,\fg]$ lies in the annhilator of $U$ and the result is just a consequence of the proposition above.\qed

The following proposition, in which we give necessary and sufficient conditions on the APK-compatible pair to yield a 2-step nilpotent Lie algebra and to assure, in such case, that the metric is flat, clarifies the conditions of the Corollary.

\begin{prop} Let $(U,H)$ be an APK-compatible pair. The Lie algebra $\fg_U$ is 2-step nilpotent if and only if $U^3=\{0\}$ and $(\mbox{\rm ann}(U))^\bot\subset \mbox{\rm ann}(U)$. If this is the case, then the pseudo-K\"ahler metric is flat if and only if $U^2\subset (U^2)^\bot$.
\end{prop}
\dem As we had seen in the proofs of Proposition \ref{flatU} and Corollary \ref{flatU2}, when $\fg_U$ is 2-step nilpotent one has $R_xR^*_y=0$ for all $x,y\in U$. But this is equivalent to $R^*_UU\subset (R^*_UU)^\bot$ and, since $R^*_UU= (\mbox{\rm ann}(U))^\bot$, it follows that $(\mbox{\rm ann}(U))^\bot\subset \mbox{\rm ann}(U)$. Recalling that we also had
$[z,[x,y]]=-(R_z+R^*_z)[x,y]$ for all $x,y,z\in U$, we arrive at 
$$0=i[z,[x,y]]+[z,[x,iy]]=-iR^*_zR^*_yx+iR^*_zR^*_xy-R^*_zR^*_{iy}x+R^*_zR^*_{x}iy=2iR^*_zR^*_xy$$
and, thus, we get $H(xzt,y)=H(t,R^*_zR^*_xy)=0$, showing that $A^3=\{0\}$. The converse follows at once from the identity $[z,[x,y]]=-(R_z+R^*_z)[x,y]$.

Since $R_xR^*_y=0$ for all $x,y\in U$, the flatness of the pseudo-K\"ahler metric is equivalent to $R_x^*R_y=0$ for all $x,y\in U$. But this is nothing but the condition
$H(U^2,U^2)=\{0\}$.\qed

\begin{example}{\em We had seen that if ${\mathcal A}$ is a nilpotent symmetric associative commutative algebra such that ${\mathcal A}^3\ne \{0\}$, then the standard pseudo-metric on $\f{aff}({\mathcal A})$ is not flat.
But one can also construct non-flat 2-step nilpotent algebras. For example, let us consider the pair $(U,H)$ where the associative algebra $U$ is the ${\mathbb C}$-span of $\{u_0,u_1,u_2,u_3,u_4\}$ with the non-trivial products
$u_4u_3=u_3u_4=u_4u_4=u_2$ and the hermitian form such that $H(u_j,u_k)=1$ if  $j+k=4$ and $H(u_j,u_k)=0$ whenever  $j+k\ne 4$. Obviously, $U^3=\{0\}$ and one has
$$ R^*_{u_3}u_2=u_0\,, \ R^*_{u_4}u_2=u_0+u_1\,, \ R^*_{u_3}u_j=R^*_{u_4}u_j=0\, ,\ \ j\neq 2 
\, .$$
Therefore, $R_xR^*_y=0$ for all $x,y\in U$. This proves that $(U,H)$ is APK-compatible and that the Lie algebra $\fg_U$ is 2-step nilpotent. However, one has
$H(U^2,U^2)\neq 0$, which shows that the metric cannot be flat.}
\end{example}

\medskip

Although there exist nilpotent Lie algebras with non flat abelian pseudo-K\"ahler structures, all of them are Ricci flat as it was shown in a more general context in \cite[Lemma 6.3]{fsp} and can be directly proved:

\smallskip

\begin{coro} Every unimodular pseudo-K\"ahler Lie algebra  with abelian complex structure is Ricci flat.
\end{coro}
\dem We had seen that for pseudo-K\"ahler Lie algebras  with abelian complex structures unimodularity and nilpotency are equivalent. Further, this means that if $(U,H)$ is an associated APK-compatible pair, then $U$ is nilpotent. This clearly implies that $R_x$ is nilpotent, hence traceless, for all $x\in U$ and, thus,  
$\mbox{Ric}(x,y)=-2\, \mbox{trace}(R_yR^*_x)=-2\,\mbox{trace}(R_{R^*_xy})=0.$\qed

\begin{example} {\rm Bearing in mind the formula given for the Ricci curvature and the explicit calculation done in the case of the algebras $\f{aff}({\mathcal A})$, one could think that Ricci flatness is only possible for nilpotent algebras. However, this is not the case. If we consider the pair $(U,H)$ where $U$ is the complex span of $\{u_0,u_1,u_2,u_3\}$ with the only non trivial products $u_1u_3=u_3u_1=u_1$, $u_3^2=u_3$ and the hermitian form defined by $H(u_0,u_3)=H(u_1,u_2)=1$ and $H(u_i,u_j)=0$ in the remaining pairs $\{i,j\}$, then one easily sees that $U$ is not nilpotent but $R_xR^*_y=R^*_xR_y=0$ for all $x,y\in U$. Hence, the pseudo-K\"ahler metric is flat. }
\end{example}

We shall finish this section con a characterization of the non Ricci flat Einstein case. Recall that a pseudo-Riemannian metric $g$ is {\it Einstein} whenever there exists $\gamma\in{\mathbb R}$ such that $g=\gamma\mbox{Ric}$.

\medskip

\begin{prop} If $(\fg,\omega,J)$ is a Lie algebra with abelian pseudo-K\"ahler structure and the pseudo-K\"ahler metric is Einstein but not Ricci flat
then $\fg = \f{aff}({\mathcal A})$, endowed with the standard abelian pseudo-K\"ahler structure, where ${\mathcal A}$ is a semisimple associative commutative
algebra and the symmetric form on ${\mathcal A}$ is (up to a scalar) the trace form of the regular representation.
\end{prop}
\dem Let us consider an APK-compatible pair $(U,H)$ such that $\fg=\fg_U$. If the metric is Einstein, then $\mbox{\rm Ric}(x,y)=-2\mbox{trace}(R_yR^*_x)$ must be non-degenerate and, hence, $\mbox{ann}(U)$ must vanish. According to Proposition \ref{annzero}, $\fg$ is holomorphically symplectomorphic to an algebra $\f{aff}({\mathcal A})$ for some symmetric associative commutative algebra $({\mathcal A},B)$. But, in this case, we had seen in the Example \ref{curv_aff} that $\mbox{Ric}((a,b),(a',b'))=-\mbox{trace}(R_{bb'})-\mbox{trace}(R_{aa'})$ for all $a,b,a',b'\in {\mathcal A}$. Bearing in mind the construction of the the pseudo-K\"ahler metric, we have for all $a,b\in {\mathcal A}$ that
$$B(a,b)=\omega ((a,0), (0,b))=g((0,-a),(0,b))=\gamma \mbox{Ric}((0,-a),(0,b))=\gamma\mbox{trace}(R_{ab}).$$
Since such form is non-degenerate if and only if ${\mathcal A}$ is semisimple, the result follows.\qed

\section{Inductive construction of pseudo-K\"ahler Lie algebras with abelian complex structure}

In \cite{m-r91} and \cite{dard}, the authors describe the method of  symplectic double extensions by one-dimensional Lie algebras to construct new symplectic Lie algebras from a given one. Actually, they prove that every nilpotent symplectic Lie algebra can be obtained by a series of double extensions starting from the zero algebra. The same constructions (slightly modified) have been used in \cite{BB-PK} to give an inductive description of all complex para-K\"ahler Lie algebras with abelian paracomplex structure. However, the method is not applicable in the pseudo-K\"ahler case unless one uses double extensions by planes instead of by lines. 

Since every pseudo-K\"ahler Lie algebra with abelian complex structure is associated to an APK-compatible pair, we will define an inductive method of construction for such pairs. The method is again based on two extensions of a given associative commutative algebra, in a similar way as it is done in \cite{baklouti} for symplectic associative commutative algebras. 

All the associative algebras of this section will be, except if the contrary is explicitly said, complex algebras. We first recall some definitions:

\begin{defi} {\em Let $U$ be an associative commutative algebra and $\varphi: U\times U\to{\mathbb C}$ a symmetric bilinear form such that $\varphi (xy,z)=\varphi (x,yz)$ for all $x,y,z\in U$. Then the vector space $U\oplus{\mathbb C}u_0$ with the product
$$(x+\alpha u_0)(x'+\alpha' u_0)=xx'+\varphi (x,x')u_0\, ,\quad x,x'\in U,\, \alpha,\alpha'\in {\mathbb C},$$
 is then an associative commutative algebra which we will call the (one-dimensional) {\it central extension of $U$ by means of $\varphi$.}}
\end{defi} 

\begin{defi} {\em Let $U$ be an associative commutative algebra, $\delta\in \f{gl}(U)$, $x_0\in U$ and let ${\mathbb C}v$ denote a one-dimensional associative algebra. On the vector space ${\mathbb C}v\oplus U$ let us define the product
$$(\alpha v+x)(\beta v+y)=\alpha\beta v^2+\alpha\beta x_0+xy+\alpha \delta y+\beta \delta x \, ,\quad \alpha,\beta\in {\mathbb C},\, x,y\in U.$$
It can be easily seen that ${\mathbb C}v\oplus U$ is an associative commutative algebra if and only if $\delta (xy)=(\delta x)y$ for all $x,y\in U$ and $\delta^2-\varepsilon \,\delta=R_{x_0}$ where $\varepsilon\in {\mathbb C}$ is such that $v^2=\varepsilon v$. In this case we will say that the algebra  ${\mathbb C}v\oplus U$ is the {\it generalized semidirect product of $U$ and ${\mathbb C}v$ by means of $(\delta,x_0)$}.
}
\end{defi} 

Our aim is to prove that, essentially, all APK-compatible pairs may be obtained by a double extension given by a central extension and a generalized semidirect product 
from another APK-compatible pair. We first give conditions to assure that a double extension of that type on an APK-compatible pair leads to another APK-compatible pair.

\begin{prop} \label{DE} Let $(U,H)$ be an APK-compatible pair, $\varepsilon\in {\mathbb C}$, $D\in \f{gl}(U)$ such that $D(xy)=(Dx)y$ for all $x,y\in U$ and $D^2-\varepsilon D=R_{a_0}$ for some $a_0\in U$ and consider a semilinear mapping $\tau:U\to U$ such that $H(x,\tau y)=\overline{H(\tau x,y)}$. Suppose that that there exist $b_0\in U$ such that  
$$\tau a_0=(D^*-\varepsilon \,\mbox{\rm id})b_0\, ,\quad \tau^2=DD^*=R_{b_0}\, ,\quad \tau b_0=Db_0$$ and that for all $x\in U$ the following conditions hold
$$ \tau R_x^*=R_x\tau\, ,\quad \tau R_x=R_x^*\tau\, ,\quad D\tau=\tau D^*\, , \quad DR_x^*=R_{\tau x}\, ,\quad \tau Dx=R^*_xb_0.$$

Let us construct a new algebra $\tilde{U}={\mathbb C}{v_0} \oplus U\oplus {\mathbb C}{ u_0}$ with a commutative product $\star$ such that $u_0$ annihilates $\tilde{U}$, $v_0\star v_0=\varepsilon v_0+a_0+\alpha_0u_0$ for some $\alpha_0\in {\mathbb C}$ and 
$$ x\star y=xy+H(x,\tau y)u_0 \, ,\quad v_0\star x=Dx+H(x,b_0)u_0 \,$$
hold for all $x,y\in U$. 
If we define a non-degenerate hermitian form $\tilde{H}$ on $\tilde{U}$ by
$$\tilde{H}(u_0,v_0)=1\, ,\quad \tilde{H}(u_0,u_0)=\tilde{H}(v_0,v_0)=\tilde{H}(u_0,x)=\tilde{H}(v_0,x)=0\, ,\quad \tilde{H}(x,y)=H(x,y)\, , \quad x,y\in U,$$
then $(\tilde{U},\tilde{H})$ is an APK-compatible pair.
\end{prop}
\dem Let us first show that $\tilde{U}$ is actually an associative commutative algebra constructed by performing a central extension  on $U$ and then a semidirect product to the extended algebra. It is straightforward to prove that the map $\varphi: U\times U\to {\mathbb C}$ given by $\varphi (x,y)=H(x,\tau y)$ is bilinear and symmetric. Moreover, for all $x,y,z\in U$ we have
$$\varphi (xy,z)=H(xy,\tau z)=H(x,R^*_y\tau z)=H(x,\tau (yz))=\varphi (x,yz).$$
This shows that $U'=U\oplus  {\mathbb C}u_0$ is the central extension of $U$ by means of $\varphi$. Now, if we consider $\delta\in {\f{gl}}(U')$ given by $\delta (x+\alpha u_0)=Dx+ H(x,b_0)u_0$ for $x\in U,\alpha\in {\mathbb C}$, then we have that
\begin{eqnarray*}
& & \delta (x)y=D(x)y+H(Dx,\tau y)u_0=D(xy)+\overline{H(\tau Dx,y)}u_0=D(xy)+\overline{H(R^*_xb_0,y)}u_0\\
& & \delta (xy+H(x,\tau y)u_0)=D(xy)+H(xy,b_0)u_0\\
& & (\delta^2-\varepsilon\, \delta)x=D^2x-\varepsilon \, Dx+H(x,D^*b_o-\varepsilon\, b_0)u_0=R_{a_0}x+H(x,\tau a_0)u_0\end{eqnarray*}
is verified for all $x,y\in U$. This shows that $\delta$ fulfills the necessary conditions to construct a generalized semidirect product of $U'$ and ${\mathbb C}v_0$. Further, one immediatly sees that the product on $\tilde{U}$ is nothing but the product on such generalized semidirect product.

It only remains to prove, thus, that the condition (\ref{compat}) of the definition of an APK-compatible pair is verified. Let us denote by $\tilde{R}_\xi$ the multiplication by $\xi\in \tilde{U}$ in $\tilde{U}$. It is obvious that $\tilde{R}^*_{u_0}=0$ since $u_0\in \mbox{ann}(\tilde{U})$ and rather long but direct calculations yield:
\begin{eqnarray*}\begin{array}{lll}
 \tilde{R}^*_{v_0}v_0=b_0+\overline{\alpha_0}u_0\, ,\quad &
 \tilde{R}^*_{v_0} {x}=D^*x+H(x,a_0)u_0 \, ,\quad & \tilde{R}^*_{v_0} {u_0}=\overline\varepsilon u_0\, ,\vspace{0.15cm}\\
\tilde{R}^*_{x}v_0=\tau x+H(b_0,x)u_0\, ,\quad &
 \tilde{R}^*_{x}y=R^*_xy+H(y,Dx)u_0\, ,\quad & \tilde{R}^*_x{u_0}=0\, ,\\\end{array}\end{eqnarray*}
for all $x,y\in U$. We then have for $x\in U$ that 
\begin{eqnarray*}
& & v_0\star \tR^*_{v_0}x=v_0\star D^*x=DD^*x+H(D^*x,b_0)u_0\\
& & x\star \tR^*_{v_0}v_0=x\star b_0=xb_0+H(x,\tau b_0)u_0
\end{eqnarray*}
and, therefore, they are equal if and only if $DD^*=R_{b_0}$ and $Db_0=\tau b_0$.
If we now take $x,y\in U$, then we get
\begin{eqnarray*}
& & x\star \tR^*_{v_0}y=x\star D^*y=xD^*y+H(x,\tau D^*y)u_0\\
& & y\star \tR^*_{v_0}x=y\star D^*x=yD^*x+H(D^*x,\tau y)u_0.
\end{eqnarray*}
From the equality $DR^*_x=R_{\tau x}$ we get $R_xD^*=R^*_{\tau x}$ and hence, bearing in mind that $\tau R_x^*=R_x\tau$, we obtain for all $z\in U$, 
$$H(xD^*y,z)=H(y,z\tau (x))=H(R^*_zy,\tau x)=H(x,\tau R^*_zy)=H(x,z\tau (y))=H(yD^*x,z).$$
Since we also have $D\tau =\tau D^*$, we deduce that  $ x\star \tR^*_{v_0}y= y\star \tR^*_{v_0}x.$ From the conditions $DR^*_x=R_{\tau x}$ and $\tau ^2=R_{b_0}$ we easily get that $x\star \tR^*_y{v_0}=v_0\star \tR^*_yx$ because
\begin{eqnarray*}
& & x\star \tR^*_y{v_0}=x\star \tau y=x\tau(y)+H(x,\tau^2y)u_0\\
& & v_0\star \tR^*_yx=v_0\star R^*_yx=DR^*_yx+H(R^*_yx,b_0)u_0.
\end{eqnarray*}
Finally, when $x,y,z\in U$ we have
$$x\star \tR^*_y{z}=x\star R^*_yz=xR^*_yz+H(x,\tau R^*_yz)u_0=zR^*_yx+H(R^*_yx,\tau z)u_0=z\star R^*_yx=z\star \tR^*_yx,$$
because $(U,H)$ is APK-compatible and $\tau R^*_y=R_y\tau$. The result now follows at once since the cases involving $u_0$ are trivial because $u_0$ and $\tR_\xi^*u_0$ are in the annihilator of $\tilde{U}$ for all $\xi\in \tilde{U}$.\qed

\begin{defi} {\em Each pair  $(\tilde{U},\tilde{H})$ constructed as in the proposition above will be called an {\it APK double extension} of the pair $(U,H)$ by a one-dimensional associative algebra. }
\end{defi} 

\medskip

\begin{lemma} Let $(U,H)$ be an APK-compatible pair. If $\mbox{\rm ann}(U)\cap ( \mbox{\rm ann}(U))^\bot\neq \{0\}$, then there exists $u_0\in\mbox{\rm ann}(U)\cap R^*_UU$, $u_0\ne 0$ and a semilinear map $\lambda: U\to {\mathbb C}$ such that
$R^*_xu_0=\lambda (x)u_0$ holds for all $x\in U$.
\end{lemma}
\dem First notice that $\mbox{\rm ann}(U)\cap R^*_UU\neq \{0\}$ since $R^*_UU=(\mbox{\rm ann}(U))^\bot$. If $z\in \mbox{\rm ann}(U)$ then, for all $x,y\in U$, we have $yR^*_xz=zR^*_xy=0$, which proves that $\mbox{\rm ann}(U)\cap R^*_UU$ is stable by all $R^*_x$, $x\in U$. Since $R^*_xR^*_y=R^*_{xy}$ and $\alpha R^*_x+\beta R^*_y=R^*_{\bar\alpha x+\bar\beta y}$ for all $\alpha,\beta\in {\mathbb C}$, $x,y\in U$, the set ${\mathcal L}=\{ {\mathcal{R}}^*_x\ :\ x\in U\}$, where ${\mathcal{R}}^*_x$ holds for the restriction of $R^*_x$ to  $\mbox{ann}(U)\cap R^*_UU$, is a complex abelian Lie subalgebra of  $\mathfrak{gl}(\mbox{ann}(U)\cap R^*_UU)$. The existence of $u_0$ and a function $\lambda : U\to {\mathbb C}$ verifying the condition is then guaranteed by Lie's theorem and the semilinearity of $\lambda$ is directly deduced from $\alpha R^*_x+\beta R^*_y=R^*_{\bar\alpha x+\bar\beta y}$.\qed

\begin{theor} Every   APK-compatible pair $(\tilde{U},\tilde {H})$ such that $\mbox{\rm ann}(\tilde{U})\cap (\mbox{\rm ann}(\tilde{U}))^\bot\neq \{0\}$ is an APK double extension of another APK-compatible pair $(U,H)$ by a one-dimensional algebra.

As a consequence, every APK-compatible pair $(\tilde{U},\tilde {H})$ is either the pair corresponding to an algebra of type $\f{aff}({\mathcal A})$ with its standard pseudo-K\"ahler structure or  can be obtained by a series of APK double extensions starting from the pair associated to an algebra $\f{aff}({\mathcal A})$. 
\end{theor}
\dem Let us apply the lemma above and consider $u_0\in\mbox{\rm ann}(\tilde{U})\cap {\tilde{R}}^*_{\tilde{U}}\tilde{U}$, $u_0\ne 0$ such that $\tilde{R}^*_\xi u_0=\lambda (\xi)u_0$  for all $\xi\in \tilde{U}$. It is clear that $\tilde{H}(u_0,u_0)=0$ since $\tR^*_{\tilde{U}}\tilde{U}=(\mbox{\rm ann}(\tilde{U}))^\bot.$ As $\tilde{H}$ is non degenerate, we may find an element $v_0\in \tilde{U}$ such that $\tilde{H}(v_0,v_0)=0$ and $\tilde{H}(u_0,v_0)=1.$ This gives a decomposition of vector spaces $\tilde{U}={\mathbb C}u_0\oplus {U}\oplus {\mathbb C}v_0$ where $U=({\mathbb C}u_0+{\mathbb C}v_0)^\bot$. The subset ${\mathcal I}=({\mathbb C}u_0)^\bot={\mathbb C}u_0\oplus {U}$ is actually an ideal of $\tilde{U}$ because for all $x\in {\mathcal I}$ and $\xi\in \tilde{U}$ one has
$$\tilde{H}(u_0,x\star \xi)=\tilde{H}(\tilde{R}_\xi^*u_0,x)=\tilde{H}(\lambda (\xi)u_0,x)=\lambda (\xi)\tilde{H}(u_0,x)=0,$$
where $\star$ denotes the product on $\tilde{U}$. Note that the restriction $H$ of $\tilde{H}$  to  ${U}\times {U}$ is non degenerate and that ${U}$ is naturally endowed with a structure of associative commutative algebra if we consider the product given by the projection to $U$ of the product in ${\mathcal I}$. Moreover, if  for all $x,y\in U$ we denote $xy\in U$ the projection of $x\star y$ to $U$, we have $x\star y=x y+\varphi(x,y)u_0$, where $\varphi$ is certain symmetric bilinear form on ${U}$ such that $\varphi(xy,z)=\varphi(x,yz)$.
Let $\tau:U\to U$ be the unique semilinear mapping such $\varphi (x,y)=H(x,\tau y)$ for all $x,y\in U$, which is guaranteed because $H$ is non-degenerate. It can be easily shown that $H(x,\tau y)=\overline{H(\tau x,y)}$ and
$\tau R_x=R_x^*\tau$ because $\varphi$ is symmetric and associative. Recall that $ {\mathcal I}={\mathbb C}u_0\oplus {U}$ is an ideal of $\tilde{U}$ and, therefore, if we denote by $D\in\f{gl}(U)$ the projection on $U$ of the restriction of $\tilde{R}_{v_0}$ to $U$, we can describe the product in $\tilde{U}$ by $u_0\in \mbox{\rm ann}(\tilde{U})$ and 
$$\begin{array}{llll} v_0\star v_0=\varepsilon v_0+a_0+\alpha_0u_0\, ,& v_0\star x=Dx+\phi (x)u_0\, ,
\quad&
x\star y=xy+H(x,\tau y)u_0\, , \quad
&  x,y\in U\end{array}$$
for some some linear form $\phi:U\to {\mathbb C}$ and certain $\varepsilon,\alpha_0\in {\mathbb C}$, $a_0\in U$. The non degeneracy of $H$ let us find $b_0\in U$ such that $\phi (x)=H(x,b_0)$ for all $x\in U$.

The associativity of $\tilde{U}$ implies that $D(xy)=(Dx)y$, $\tau Dx=R^*_xb_0$ hold for all $x,y\in U$ and that $\tau a_0=D^*b_0-\varepsilon b_0$, $D^2-\varepsilon\, D=R_{a_0}$ and, since the product $\star$ is given as the one in Proposition \ref{DE}, one sees by simply reversing the arguments of its proof that the condition $\xi\tilde{R}^*_\nu\xi'=\xi'\tilde{R}^*_\nu\xi$ for $\xi,\xi',\nu\in \tilde{U}$ implies that $(U,H)$ is an APK-compatible pair and, further, that all the conditions on $\tau, D, a_0$ and $b_0$ to have an APK double extension are fulfilled.

The second part of the statement follows by applying sucessively the first part. If $(\tilde{U},\tilde {H})$ is not the pair corresponding to an algebra $\f{aff}({\mathcal A})$ then, according to Proposition \ref{aff_ext}, $\mbox{\rm ann}(\tilde{U})\cap (\mbox{\rm ann}(\tilde{U}))^\bot\neq \{0\}$ and it is an APK double extension of certain pair $(U_1,H_1)$. The same argument can be then applied to $U_1$ and so on.\qed

\begin{rema}{\em The following facts are remarkable:
\begin{enumerate}
\item[(1)] The process to view an APK-compatible pair as a series of double extensions stops when we arrive at an algebra in which the annhilator does not intersect its orthogonal. This means that the symmetric real associative commutative algebra ${\mathcal A}$ of the second part of the theorem can be taken to be one of the algebras considered in Proposition \ref{aff_ext}. In particular, every nilpotent APK-compatible pair $(U,H)$ such that $U^2\ne\{0\}$ is obtained by a sequence of APK double extensions starting from an algebra with zero multiplication.
\item[(2)] When $\fg_U$ is the Lie algebra constructed with an APK-compatible pair $(U,H)$, the signature of the pseudo-K\"ahler metric is twice the one of $H$ and when $(\tilde{U},\tilde{H})$ is an APK double extension of a compatible pair $(U,H)$ then the signature of $\tilde H$ is given by $\mbox{sig}(\tilde H)=\mbox{sig} (H)+1$. This shows that if the pseudo-K\"ahler metric is definite positive, then $\fg$ cannot be an APK double extension and must actually be an algebra  $\f{aff}({\mathcal A})$ for an associative commutative symmetric algebra $({\mathcal A},B)$ as those in Proposition \ref{aff_ext}, this is to say,  an orthogonal sum ${\mathcal A}=\mbox{\rm ann}({\mathcal A})\oplus {\mathcal A}_1\oplus\cdots \oplus {\mathcal A}_p$ where each ideal ${\mathcal A}_i$ is of the form  $e_i{\mathcal A}$ for some idempotent $e_i$. When the K\"ahler metric of $\f{aff}({\mathcal A})$ is definite positive,  also $B$ must be definite (negative). If the nilradical ${\mathcal N}_i$ of ${\mathcal A}_i$ is non-zero, $x\in {\mathcal N}_i$ and $k\in {\mathbb N},k\ge 2$ is such that $x^{2k}=0,x^{2k-2}\ne 0$ then we have $0=B(x^{2k},x^2)=B(x^{k+1},x^{k+1})$ and hence $x^{k+1}=0$, which is obviously a contradiction with $x^{2k-2}\ne 0$ if $k\ge 3$ and also if $k=2$ because $x^3=0$ implies $B(x^2,x^2)=B(x^3,x)=0$. Thus ${\mathcal N}_i=\{0\}$ and then $e_i{\mathcal A}$ must be simple. This proves that each ideal $e_i{\mathcal A}$ is isomorphic to either ${\mathbb R}$ or ${\mathbb C}$. But a simple calculation shows that ${\mathbb C}$ does not admit a definite form of the type defined in Definition \ref{def_symm}. Thus, one immediately obtains the well-known result (see, for instance, \cite[Theorem 4.1]{Andrada3}) that a K\"ahler Lie algebra with abelian complex structre is an orthogonal sum of several copies of $\f{aff}({\mathbb R})$ and an even dimensional abelian algebra.
\end{enumerate}}\end{rema}

\medskip

\noindent {\it Acknowledgments:} We would like to thank Professor Sa\"{\i}d Benayadi for his multiple and interesting suggestions.

\end{document}